\documentclass{amsart}
\usepackage{amscd}
\usepackage{amssymb}
\newtheorem{theorem}{Theorem}
\numberwithin{equation}{subsection}
\DeclareMathOperator{\ssc}{sc}
\DeclareMathOperator{\val}{val}

\begin{document}

\title{Regular points in affine Springer fibers}
\date{August 7, 2003}

\author{Mark Goresky}
\address{M. Goresky, School of Mathematics, IAS, Princeton, NJ 08540}
\email{goresky@ias.edu}
\thanks{M.G. was partially supported by NSF Grant DMS-0139986}

\author{Robert Kottwitz}
\address{R. Kottwitz, Department of Mathematics, University of Chicago, 5734
University Ave., Chicago, IL 60637}
\email{kottwitz@math.uchicago.edu}
\thanks{R.K. was partially supported by NSF Grant DMS-0245639}

\author{Robert MacPherson}
\address{R. MacPherson, School of Mathematics, IAS, Princeton, NJ 08540}
\email{rdm@ias.edu}



\maketitle

\section{Introduction}
Let $G$ be a connected reductive group over $\mathbf C$ with Lie algebra
$\mathfrak g$. We put $F=\mathbf C((\epsilon))$ and $\mathcal
O=\mathbf C[[\epsilon]]$. Let $X=X_G$ denote the affine Grassmannian
$G(F)/G(\mathcal O)$.  For $u \in \mathfrak g(F)$ we write $X^u$ for the
affine Springer fiber 
\begin{equation*} X^u=\{g \in G(F)/G(\mathcal O) : \mathrm{Ad}(g^{-1})(u)
\in \mathfrak g(\mathcal O)\}.
\end{equation*}
studied by Kazhdan and Lusztig in \cite{kazhdan-lusztig88}.

For $x=gG(\mathcal O) \in X^u$ the $G(\mathcal O)$-orbit (for the adjoint
action) of $\mathrm{Ad}(g^{-1})(u)$ in $\mathfrak g (\mathcal O)$ depends
only on~$x$, and its image under $\mathfrak g(\mathcal O) \twoheadrightarrow
\mathfrak g (\mathbf C)$ is a well-defined $G(\mathbf C)$-orbit in
$\mathfrak g(\mathbf C)$. We say that $x \in X^u$ is \emph{regular} if the
associated orbit is regular in $\mathfrak g(\mathbf C)$. (Recall that an
element of~$\mathfrak g(\mathbf C)$ is regular if the nilpotent part of its
Jordan decomposition is a principal nilpotent element in the centralizer of
the semisimple part of its Jordan decomposition.) We write $X^u_{\mathrm
{reg}}$ for the (Zariski open) subset of regular elements in $X^u$. 

{}From now on we assume that $u$ is regular semisimple with
centralizer~$T$, a maximal torus in~$G$ over~$F$. Assume further that $u$
is \emph{integral}, by which we mean that $X^u$ is non-empty. Kazhdan and
Lusztig
\cite{kazhdan-lusztig88} show that $X^u$ is then a locally finite union of
projective algebraic varieties, and in Cor.~1 of \S4 of
\cite{kazhdan-lusztig88} they show that the open subset $X^u_{\mathrm
{reg}}$ of~$X^u$ is non-empty (and hence dense in at least one irreducible
component of~$X^u$). 
The action of~$T(F)$ on~$X$ clearly preserves the subsets $X^u$ and 
$X^u_{\mathrm {reg}}$.  Bezrukavnikov \cite{bezrukavnikov96} proved
that $X^u_{\mathrm {reg}}$ forms a single orbit under $T(F)$. (Actually
Kazhdan-Lusztig and Bezrukavnikov consider only topologically nilpotent
elements $u$, but the general case can be reduced to their
special case by using the topological Jordan decomposition of~$u$.)

The goal of this paper is to characterize regular elements in~$X^u$ (for
integral regular semisimple $u$ as above). When $T$ is elliptic (in other
words, $F$-anisotropic modulo the center of~$G$) the characterization gives
no new information. At the other extreme, in the split case, the
characterization gives a clear picture of what it means for a point in
$X^u$ to be regular.  

We will now state our characterization in the split case, leaving the more
technical general statement to the next section (see Theorem
\ref{mainthm}). Fix a split maximal torus $A \subset G$ over~$\mathbf C$
and denote by $\mathfrak a$ its Lie algebra. We identify the affine
Grassmannian $A(F)/A(\mathcal O)$ for~$A$ with the cocharacter lattice
$X_*(A)$, the cocharacter $\mu$ corresponding to the class of
$\mu(\epsilon)$ in $A(F)/A(\mathcal O)$. For any Borel subgroup $B=AN$
containing~$A$ ($N$ denoting the unipotent radical of~$B$) there is a
well-known retraction $r_B:X \to X_*(A)$ defined using the Iwasawa
decomposition $G(F)=N(F)A(F)G(\mathcal O)$: the fiber of~$r_B$ over $ \mu
\in X_*(A)$ is  $N(F)\mu(\epsilon)G(\mathcal O)/G(\mathcal
O)$. The family of cocharacters $r_B(x)$ ($B$ ranging through all Borel
subgroups containing~$A$) has been studied by Arthur \cite[Lemma
3.6]{arthur76}; it is the volume of the convex hull of these points that
arises as the weight factor for (fully) weighted orbital integrals for
elements in $A(F)$. In particular Arthur shows that for $x \in X$ and any
pair $B$, $B'$ of adjacent Borel subgroups containing $A$, there is a unique
non-negative integer $n(x,B,B')$ such that
\begin{equation}
r_B(x)-r_{B'}(x)=n(x,B,B') \cdot \alpha^\vee_{B,B'}
\end{equation}
where $\alpha_{B,B'}$ is the unique root of~$A$ that is positive for~$B$ and
negative for $B'$. 

The main result of this paper (in the split case) is that for $x \in X^u$
\begin{equation}\label{ineq}
n(x,B,B') \le \val \alpha_{B,B'}(u)
\end{equation}
for every pair $B$, $B'$ of adjacent Borel subgroups containing~$A$, and
that $x \in X^u$ is regular if and only if all the inequalities \eqref{ineq}
are actually equalities.

\section{Statements}
\subsection{Notation}  We write $\mathfrak g$ for the Lie algebra of $G$ and 
follow the same convention for groups denoted by other letters. 
 
Choose an algebraic closure 
$\overline F$  of $F$ and let $\Gamma=\mathrm{Gal}(\overline F/F)$. We write
$G_F$ for the $F$-group obtained from $G$ by extension of scalars from
$\mathbf C$ to $F$. 

As before we use $\mu \mapsto \mu(\epsilon)$ to identify the
cocharacter group $X_*(A)$ with $A(F)/A(\mathcal O)$. By means of this
identification the canonical surjection $A(F) \to A(F)/A(\mathcal O)$ can be
viewed as a surjection
\begin{equation}\label{map.A}
A(F) \to X_*(A).
\end{equation}
Let $\Lambda=\Lambda_G$ denote the
quotient of the coweight lattice
$X_*(A)$ by the coroot lattice (the subgroup of $X_*(A)$ generated by the
coroots of $A$ in $G$). Up to canonical isomorphism $\Lambda$ is independent
of the choice of $A$; moreover when defining $\Lambda$ we could replace $A$
by any maximal torus $T$ in $G_F$. There is a canonical surjective 
homomorphism
\begin{equation}\label{map.G}
G(F) \to \Lambda,
\end{equation}
characterized by the following two properties: it is trivial on the image of
$G_{\ssc}(F)$ in $G(F)$ (where $G_{\ssc}$ denotes the simply connected cover
of the derived group of~$G$), and its restriction to $A(F)$ coincides with
the composition of \eqref{map.A} and the canonical surjection  $X_*(A)
\to \Lambda$. 
 
Recall that
$X$ denotes the affine Grassmannian 
$G(F)/G(\mathcal O)$ for
$G$. The homomorphism \eqref{map.G} is trivial on~$G(\mathcal O)$ and hence
induces a canonical surjection  
\begin{equation}
\nu_G:X \to \Lambda,
\end{equation}
whose fibers are the connected components of~$X$. 

\subsection{Parabolic subgroups} We will concerned with parabolic subgroups
$P$ of~$G$ containing $A$. Such a parabolic subgroup has a unique Levi
subgroup~$M$ containing $A$, and we refer to $M$ as \emph{the} Levi
component of~$P$. 

As usual, by a Levi subgroup of $G$, we mean a Levi subgroup of some
parabolic subgroup of~$G$. Let
$M$ be a Levi subgroup of
$G$  containing
$A$. We write $\mathcal P(M)$ for the set of parabolic subgroups of $G$ that
contain $A$ and have Levi component~$M$. Thus any $P \in \mathcal P(M)$
can be written as $P=MN$ where
$N=N_P$ denotes the unipotent radical of $P$. As usual there is a notion of
adjacency: two parabolic subgroups $P=MN$ and $P'=MN'$ in $\mathcal P(M)$
are said to be \emph{adjacent} if there exists (a unique) parabolic subgroup
$Q=LU$ containing both $P$ and $P'$ such that the semisimple rank of $L$ is
one greater than the semisimple rank of $M$. Thus $U=N\cap N'$, and, 
moreover, if $L$ is chosen so that $L \supset A$, then 
\begin{equation*}
\mathfrak l=\mathfrak m\oplus (\mathfrak n \cap \bar{\mathfrak n'}) \oplus
(\mathfrak n' \cap \bar {\mathfrak n})
\end{equation*} 
where $\bar N$ denotes the unipotent radical of the
parabolic subgroup
$\bar P=M\bar N$ opposite to $P$ (and where  $\bar{ N}'$ is opposite to
$N'$).

Given adjacent $P$, $P'$ in $\mathcal P(M)$ we define an element
$\beta_{P,P'} \in \Lambda_M$ (the coweight lattice for $A$ modulo the coroot
lattice for $M$) as follows. Consider the collection of elements in
$\Lambda_M$ obtained from coroots $\alpha^\vee$ where $\alpha$ ranges
through the set of roots of
$A$ in $\mathfrak n \cap \bar{\mathfrak n}'$. We define $\beta_{P,P'}$ to be
the unique element in this collection such that all other members in the
collection are positive integral multiples of $\beta_{P,P'}$. Note that
although $\Lambda_M$ may have torsion elements, the elements in our
collection lie in the kernel of the canonical map from $\Lambda_M$ to
$\Lambda_G$, and this kernel is torsion-free. Thus any member of our
collection can be written uniquely as a positive integer times
$\beta_{P,P'}$. Note also that $\beta_{P',P}=-\beta_{P,P'}$. In case $M=A$,
so that $P,P'$ are Borel subgroups, $\beta_{P,P'}$ is the unique  coroot of
$A$ that is positive for
$P$ and  negative for $P'$. 

\subsection{Retractions from $X$ to $X_M$}\label{retractions}
 The inclusion of $M(F)$ into $G(F)$ induces an inclusion of the affine
Grassmannian $X_M$ for $M$ into the affine Grassmannian $X$ for $G$. Let $P
\in
\mathcal P(M)$ and  let $X_P$ denote the set $P(F)/P(\mathcal O)$. The
canonical inclusion of
$P$ in $G$ induces a bijection $ i$ from $X_P$ to $X$, and the canonical
surjection $P \to M$ induces a canonical surjective map $p$ (of sets) from
$X_P$ to $X_M$. We define the retraction
$r_P=r^G_P:X \to X_M$ as the composed map $p \circ i^{-1}$. Given $x \in X$
we often denote by
$x_P$ the image of $x$  under the retraction $r_P$.

These retractions satisfy the following transitivity property. Suppose that
$L\supset M$ are Levi subgroups containing $A$, and suppose further that $P
\in \mathcal P(M)$ and $Q \in \mathcal P(L)$ satisfy
$Q \supset P$. Let $P_L$ denote the parabolic subgroup $P\cap L$ in $L$.
Then 
\begin{equation} r^G_P=r^L_{P_L}\circ r^G_Q.
\end{equation} Moreover, for any  $x \in X$ the element $\nu_M(x_P)$ maps to
$\nu_L(x_Q)$ under the canonical surjection $\Lambda_M \to \Lambda_L$, and
in particular 
$\nu_M(x_P) \mapsto \nu_G(x)$ under  $\Lambda_M \to \Lambda_G$. 
 
\subsection{Definition of $n(x,P,P')$} A point $x \in X$ determines  points
$\nu_M(x_P)$ in $\Lambda_M$, one for each $P \in \mathcal P(M)$. This family
of points arises in the definition of the weighted orbital integrals
occurring in Arthur's work. A basic fact \cite{arthur76} about this family of
points is that whenever $P$, $P'$ are adjacent parabolic subgroups in
$\mathcal P(M)$, there is a (unique) non-negative integer $n(x,P,P')$ such
that
\begin{equation}
\nu_M(x_P)-\nu_M(x_{P'})=n(x,P,P')\cdot \beta_{P,P'}.
\end{equation}  

\subsection{Fixed point sets $X^u$} Let $u \in \mathfrak g(F)$.  Define
a subset
$X^u$ of $X$ by
\begin{equation*} X^u=\{g \in G(F)/G(\mathcal O) : \mathrm{Ad}(g^{-1})(u)
\in \mathfrak g(\mathcal O)\}.
\end{equation*}

\subsection{Conjugacy classes associated to fixed points}\label{conjugacy} 
Let $u \in \mathfrak g(F)$.
 Suppose that the coset $x=gG(\mathcal O)$ lies in $X^u$. The image of
$\mathrm{Ad}(g^{-1})(u)$ under the canonical surjection $\mathfrak
g(\mathcal O) \to
\mathfrak g(\mathbf C)$ gives a well-defined $G(\mathbf C)$-conjugacy class
$\bar u_G(x)$ (for the adjoint action) in $\mathfrak g(\mathbf C)$. 

As above let $M$ be  a Levi subgroup of $G$ and let $P \in \mathcal P(M)$.
Now suppose that $u
\in \mathfrak m(F)$ and that $x \in X^u$. Choose $p \in P(F)$ such that
$x=pG(\mathcal O)$; thus $x_P$ is the coset $m M(\mathcal O)$, where $m$
denotes the image of $p$ under the canonical homomorphism from $P$ onto $M$.
Of course $\mathrm{Ad}(p^{-1})(u)$ lies in $\mathfrak p(\mathcal O)$, and
its image  in $\mathfrak p(\mathbf C)$ gives a well-defined $P(\mathbf
C)$-conjugacy class  
$\bar u_P(x)$ in $\mathfrak p(\mathbf C)$. It follows that $x_P$ lies in
$X^u_M$ (as was first noted by Kazhdan-Lusztig \cite{kazhdan-lusztig88}), and
also that
$\bar u_P(x)$ maps to $\bar u_G(x)$ (respectively,
$\bar u_M(x_P)$) under the map on conjugacy classes induced by  $\mathfrak
p(\mathbf C) \hookrightarrow
\mathfrak g(\mathbf C)$ (respectively, $\mathfrak p(\mathbf C)
\twoheadrightarrow
\mathfrak m(\mathbf C)$).  

\subsection{Review of regular elements} An element $u \in \mathfrak
g(\mathbf C)$ is
\emph{regular} if   the nilpotent part of its Jordan decomposition is a
principal nilpotent element in the centralizer of the semisimple part of its
Jordan decomposition, or, equivalently, if the set of Borel subalgebras
containing $u$ is finite.  It is well-known that the set of regular elements
in 
$\mathfrak g(\mathbf C)$ is open. 

As above let $M$ be  a Levi subgroup of $G$ and let $P \in \mathcal P(M)$.
Suppose that $u$ is a regular element in $\mathfrak g(\mathbf C)$ that
happens to lie in $\mathfrak p(\mathbf C)$. Then the image $u_M$ of $u$ in
$\mathfrak m(\mathbf C)$ is regular in $\mathfrak m(\mathbf C)$. 
\subsection{Regular points in $X^u$} We say that $x \in X^u$ is
\emph{regular} if the associated conjugacy class $\bar u_G(x) \in \mathfrak
g(\mathbf C)$ consists of regular elements.  We denote by
$X^u_{\mathrm{reg}}$ the set of regular elements in $X^u$; the subset
$X^u_{\mathrm {reg}}$ is open in $X^u$.

As above let $M$ be  a Levi subgroup of $G$ and let $P \in \mathcal P(M)$.
Suppose that $u \in
\mathfrak  m(F)$. We have already seen that $r_P$ maps $X^u$ into $X^u_M$,
and that the conjugacy class in
$\mathfrak g(\mathbf C)$ associated to 
$x \in X^u$ is compatible with the conjugacy class in $\mathfrak m(\mathbf
C)$ associated to 
 the retracted  point $x_P \in X^u_M$, compatible in the sense that there is
a conjugacy class in
$\mathfrak p(\mathbf C)$ that maps to both of them. Therefore $x_P$ is
regular in $X^u_M$ if $x$ is regular in
$X^u$.

\subsection{Set-up for the main result}\label{set-up} As before let $M$
denote a Levi subgroup of
$G$ containing
$A$. We now assume that $u$ is an integral regular semisimple element of
$\mathfrak g(F)$ that happens to lie in $\mathfrak m (F)$. (It is equivalent
to assume that the centralizer
$T$ of
$u$ is contained in $M_F$.) For each pair $P=MN$, $P'=MN'$ of adjacent
parabolic subgroups in
$\mathcal P(M)$ we are going to define a non-negative integer $n(u,P,P')$.
This collection of integers measures how far $X^u$ sticks out from $X_M^u$. 

As before we  need the parabolic subgroups $\bar P=M\bar N$ and  $\bar
P'=M\bar N'$ opposite to $P$ and $P'$ respectively. Let $\alpha$ be a  root
of $T$ in $N \cap \bar N'$. Since $T$, $N$ and
$N'$ are defined over $F$, the group $\mathrm{Gal}(\overline F/F)$ preserves
the set of roots  of $T$ in $N \cap \bar N'$. Let
$F_{\alpha}$ denote the field of definition of $\alpha$, so that 
$\mathrm{Gal}(\overline F/F_\alpha)$ is the stabilizer of $\alpha$ in 
$\mathrm{Gal}(\overline F/F)$. For any finite extension $F'$ of $F$
(\emph{e.g.} $F_\alpha$) we normalize the valuation
$\mathrm{val}{_{F'}}$ on
$F'$ so that a uniformizing element in
$F'$ has valuation 1, or, equivalently, so that $\epsilon$ has valuation
$[F':F]$. There exists a unique positive integer $m_\alpha$ such that the
image of the element
$\alpha^\vee$ in $\Lambda_M$ is equal to $m_\alpha \cdot
\beta_{P,P'}$, where $\beta_{P,P'}$ is the element of $\Lambda_M$ defined
above. Note that
$m_\alpha$ depends only on the orbit of $\alpha$ under the Galois group;
here we use that the Galois group acts on the cocharacter group of $T$
through the Weyl group of $M$, so that any two elements in the Galois orbit
of $\alpha^\vee$ have the same image in $\Lambda_M$.  
 Finally we define $n(u,P,P')$ as the sum
\begin{equation} n(u,P,P')=\sum \mathrm{val}_{F_\alpha}(\alpha(u)) \cdot
m_\alpha,
\end{equation} where  the sum is taken over a set of representatives $\alpha$
of the orbits of
$\mathrm{Gal}(\overline F/F)$ on the set of roots  of $T$ in $N \cap \bar
N'$. In the special case that $M=A$ (and hence $T=A$) $n(u,P,P')$ is equal
to $\mathrm{val}_F(\alpha(u))$, where 
$\alpha$ is the unique root of $A$ that is positive for $P$ and negative for
$P'$.  
\begin{theorem}\label{mainthm}
Let $M$ and $u$ be as above, and let $x \in X^u$. Recall that
$x_P \in X_M^u$ for all
$P
\in \mathcal P(M)$. 
\begin{itemize}
\item[(a)] For every pair $P,P' \in\mathcal P(M)$ of adjacent parabolic
subgroups  
 $$n(x,P,P') \le n(u,P,P').$$

\item[(b)] The point $x$ is regular in $X^u$ if and only if the following
two conditions hold: 
\itemize
\item[(i)] the point $x_P$ is regular in $X^u_M$ for all $P \in
\mathcal P(M)$, and 
\item[(ii)] for every pair $P,P'\in \mathcal P(M)$ of adjacent parabolic
subgroups   $$n(x,P,P')  =  n(u,P,P').$$
\end{itemize}
\end{theorem}

\section{Proofs}
\subsection{The case of $SL(2)$}\label{SL(2)} The key step in proving our
main theorem is to verify it for $SL(2)$, where it reduces to a computation
that can be found in \cite{langlands80}. To keep
things self-contained we reproduce the calculation here. Let $A$, $B$, $\bar
B$ denote the diagonal, upper triangular and lower triangular subgroups of
$SL(2)$ respectively, and let $\alpha$ be the unique root of $A$ that is
positive for $B$. Of course $\beta_{B,\bar B}=\alpha^\vee$.   Let $x \in X$
and let $u = \begin{bmatrix} c & 0 \\ 0 & -c \end{bmatrix}$ for non-zero $c
\in
\mathcal O$. Note that $n(u,B,\bar B)=\mathrm{val}_F(c)$. We will show that
$x \in X^u$ if and only if
$n(x,B,\bar B)\le n(u,B,\bar B)$, and that $x
\in X^u_{\mathrm{reg}}$ if and only if $n(x,B,\bar B) = n(u,B,\bar B)$.

The difference $\nu_A(x_B)-\nu_A(x_{\bar B})$ and the sets $X^u$ and
$X^u_{\mathrm{reg}}$ are invariant under the action of $A(F)$ on $X$, so it
is enough to consider $x$ of the form $x=gG(\mathcal O)$ with 
$g=\begin{bmatrix}1&0\\t&1\end{bmatrix}$. (Note that for
this reason our calculations  apply just as well to any group whose
semisimple rank is $1$.) For such
$x$ we have
$\nu_A(x_{\bar B})=0$. If $t \in \mathcal O$, then $\nu_A(x_{B})=0$. If $t
\notin \mathcal O$, then
$\begin{bmatrix} 0&-1\\1&t^{-1}\end{bmatrix} \in G(\mathcal O)$ and thus
\begin{equation*}
\begin{bmatrix} 1&0\\t&1
\end{bmatrix}=
\begin{bmatrix}t^{-1}&1\\0&t
\end{bmatrix}
\begin{bmatrix}0&-1\\1&t^{-1}
\end{bmatrix}\in
\begin{bmatrix}t^{-1}&1\\0&t
\end{bmatrix}\cdot G(\mathcal O),
\end{equation*} which shows that $\nu_A(x_B)=\mathrm{val}_F(t^{-1})\cdot
\alpha^\vee$. We conclude that
$n(x,B,\bar B)$ equals $0$ if $t \in \mathcal O$ and equals
$\mathrm{val}_F(t^{-1})$ if $t \notin
\mathcal O$. In any case $n(x,B,\bar B)$ is a non-negative integer. 

For $x,u$ as above we have
\begin{equation*} \mathrm{Ad}(g^{-1})u=
\begin{bmatrix}c&0\\-2ct&-c
\end{bmatrix}.
\end{equation*} Therefore $x \in X^u$ $\iff$ $ct \in \mathcal O$ $\iff$
$n(x,B,\bar B) \le n(u,B,\bar B)$. Moreover 
$x \in X^u_{\mathrm{reg}}$ $\iff$ $ct \in \mathcal O^\times$ or  ($c \in
\mathcal O^\times$ and $t
\in \mathcal O$) $\iff$
$n(x,B,\bar B) = n(u,B,\bar B)$.

\subsection{Review of $n(x,P,P')$}\label{integrality}
 We need to review Arthur's proof of the existence of the non-negative
integers
$n(x,P,P')$. We
begin with the case
$M=A$. Let
$x \in X$. We must check that for any two adjacent Borel subgroups $P$, $P'
\in \mathcal P(A)$ there is a (unique) non-negative integer $n(x,P,P')$ such
that
\begin{equation*}
\nu_A(x_P)-\nu_A(x_{P'})=n(x,P,P')\cdot \alpha^\vee,
\end{equation*} where $\alpha$ is the unique root of $A$ that is positive
for $P$ and negative for $P'$. For this we consider the unique parabolic
subgroup $Q$ containing $P$ and $P'$ whose Levi component $L$ has semisimple
rank $1$. By transitivity of retractions we have
\begin{equation}\label{2.2.1}
\nu_A(x_P)-\nu_A(x_{P'})=\nu_A(y_B)-\nu_A(y_{\bar B})
\end{equation} where $y=x_Q$ and $B=L\cap P$, $\bar B=L \cap P'$. This
reduces us to the case in which $G$ has semisimple rank $1$, which has
already been done. For future use we note that \eqref{2.2.1} can be
reformulated as the equality
\begin{equation*} n(x,P,P')=n(y,B,\bar B).
\end{equation*}

Again let $x \in X$. Now we check that for any Levi subgroup $M \supset A$
and any adjacent parabolic subgroups $P=MN$,
$P'=MN'$ in $\mathcal P(M)$ there is a (unique) non-negative integer
$n(x,P,P')$ such that
\begin{equation*}
\nu_M(x_P)-\nu_M(x_{P'})=n(x,P,P')\cdot \beta_{P,P'}.
\end{equation*} Fix a Borel subgroup $B_M$ in $M$ and let $B$ (respectively,
$B'$) be the inverse image of $B_M$ under $P \to M$ (respectively, $P' \to
M$); thus $B$, $B'$ are Borel subgroups containing $A$. 

Now choose a minimal gallery of Borel subgroups $B=B_0,B_1,B_2,\dots,B_l=B'$
joining $B$ to $B'$, and for $i=1,\dots,l$ let $\alpha_i$ be the unique root
of $A$ that is positive for $B_{i-1}$ and negative for $B_i$. Then 
\begin{equation*}
\nu_A(x_B)-\nu_A(x_{B'})=\sum_{i=1}^l n(x,B_{i-1},B_i)\cdot \alpha_i^\vee.
\end{equation*}
 Note that $\{\alpha_1,\dots,\alpha_l\}$ is precisely the set of roots of
$A$ in $\mathfrak n \cap
\bar{\mathfrak n}'$ and that for each $i$ there exists a (unique) positive
integer $m_i$ such that the image of $\alpha_i^\vee$ in $\Lambda_M$ is equal
to $m_i \cdot \beta_{P,P'}$. Applying the canonical surjection $\Lambda_A
\to \Lambda_M$ to the previous equation, we find that (see
\ref{retractions}) 
\begin{equation*}
\nu_M(x_P)-\nu_M(x_{P'})=n(x,P,P')\cdot \beta_{P,P'}.
\end{equation*} where $n(x,P,P')$ is the non-negative integer
\begin{equation*}
\sum_{i=1}^l m_i \cdot n(x,B_{i-1},B_i).
\end{equation*}

\subsection{Proof of part of the main theorem in case $A=T$} Let $u \in
\mathfrak a(\mathcal O)$ and assume that $u$ is regular in $\mathfrak g(
F)$.  Let $x \in X^u$. 

Let $M$ be a Levi subgroup of $G$ containing $A$. We are now going to prove
the first main assertion in our theorem, namely that for any pair of
adjacent $P,P' \in \mathcal P(M)$ there is an inequality
\begin{equation*} n(x,P,P') \le n(u,P,P').
\end{equation*}

Let $B,B',B_0,\dots,B_l$ and $\alpha_i$, $m_i$ $(i=1,\dots,l$) be as in
\ref{integrality}. Then by definition
\begin{equation*} n(u,P,P') =\sum_{i=1}^l m_i \cdot
\mathrm{val}_F(\alpha_i(u)).
\end{equation*}

Let $M_i$ be the Levi subgroup containing $A$ whose root system is $\{\pm
\alpha_i\}$, and let
$B'_{i-1}$, $B'_i$ denote the Borel subgroups in $M_i$ obtained by
intersecting $B_{i-1}$, $B_i$ with $M_i$. Let $Q_i$ be the unique parabolic
subgroup in $\mathcal P(M_i)$ such that $Q_i$ contains $B_{i-1}$ and $B_i$.
We showed in
\ref{integrality} that 
\begin{equation*} n(x,P,P') =\sum_{i=1}^l m_i \cdot n(x,B_{i-1},B_i)
\end{equation*} and that
\begin{equation*}
 n(x,B_{i-1},B_i)=n(y_i,B'_{i-1},B'_i),
\end{equation*} where $y_i=x_{Q_i} \in X^u_{M_i}$. Since $M_i$ has
semisimple rank 1, we know that 
\begin{equation*} n(y_i,B'_{i-1},B'_i)\le \mathrm{val}_F(\alpha_i(u)).
\end{equation*} This completes the proof of the first main assertion. 

Now suppose that $x$ is regular in $X^u$. Then each point $y_i \in
X^u_{M_i}$ above is regular in
$X^u_{M_i}$, and therefore from the rank $1$ case (see \ref{SL(2)}) we know
that 
\begin{equation*} n(y_i,B'_{i-1},B'_i) = \mathrm{val}_F(\alpha_i(u)).
\end{equation*} We conclude that if $x $ is regular in $X^u$, then 
\begin{equation*} n(x,P,P') = n(u,P,P'),
\end{equation*} which is another of the assertions in our theorem. 

\subsection{Proof of the rest of the main theorem in case $M=A=T$}

We continue with  $u \in \mathfrak a(\mathcal O)$ and $x \in X^u$ as before,
but for the moment we only consider the case $M=A$.  We assume that
\begin{equation}\label{*} n(x,P,P')=\mathrm{val}_F(\alpha_{P,P'}(u)) 
\end{equation} for all adjacent Borel subgroups $P,P' \in \mathcal P(A)$,
where $\alpha_{P,P'}$ is the unique root of $A$ that is positive for $P$ and
negative for
$P'$. We want to prove that $x$ is regular in $X^u$. To do so we must first
select a suitable Borel subgroup $B \in \mathcal P(A)$. 

Let $u_0 \in \mathfrak a(\mathbf C)$ denote the image of $u$ under
$\mathfrak a(\mathcal O) \to 
\mathfrak a(\mathbf C)$, and let $M$ denote the centralizer of $u_0$ in $G$.
Thus $M$ is a Levi subgroup of $G$ containing $A$, and we choose $P \in
\mathcal P(M)$. Then we obtain a suitable Borel subgroup by taking any $B
\in \mathcal P(A)$ such that $B \subset P$.   For any $B$-simple root
$\alpha$ we denote by $B_\alpha$ the unique Borel subgroup in $\mathcal
P(A)$ that is adjacent to $B$ and for which $\alpha$ is negative, and we
write $P_\alpha$ for the unique parabolic subgroup containing $B$ and
$B_\alpha$ such that the semisimple rank of the Levi component $M_\alpha$ of
$P_\alpha$ is $1$.   Consider the element (well-defined up to $B(\mathbf
C)$-conjugacy) $v:=\bar u_B(x) \in \mathfrak b(\mathbf C)$ defined in
\ref{conjugacy}. The equation \eqref{*} plus the semisimple rank $1$ theory
implies that the points $x_{P_\alpha} \in X^u_{M_\alpha}$ are regular, and
this in turn implies  (see \ref{conjugacy}) that for every $B$-simple root
$\alpha$ the image of the element
$v$ under
$\mathfrak b(\mathbf C) \hookrightarrow \mathfrak p_\alpha(\mathbf C)
\twoheadrightarrow
\mathfrak m_\alpha(\mathbf C)$ is regular in $\mathfrak m_\alpha(\mathbf
C)$. Moreover it is evident that the image of
$v$ under the canonical surjection $\mathfrak b(\mathbf C)
\twoheadrightarrow \mathfrak a(\mathbf C)$ is equal to $u_0$. Using only
these facts, we now check that $v$ is regular in
$\mathfrak g(\mathbf C)$ (and hence that $x$ is regular in $X^u$). 

Let $v=v_s+v_n$ be the Jordan decomposition of $v$, with $v_s$ semisimple
and $v_n$ nilpotent. Since it is harmless to replace $v$ by any $B(\mathbf
C)$-conjugate, we may assume without loss of generality that $v_s \in
\mathfrak a(\mathbf C)$. Then, since
$v_s\mapsto u_0$ under
$\mathfrak b(\mathbf C)
\twoheadrightarrow
\mathfrak a(\mathbf C)$, it follows that $v_s=u_0$. Since $v_n$ commutes
with $v_s=u_0$, it lies in
$\mathfrak m(\mathbf C)$, and we must check that $v_n$ is a principal
nilpotent element in 
$\mathfrak m(\mathbf C)$. As
$v_n$ lies in the Borel subalgebra
$(\mathfrak b
\cap
\mathfrak m)(\mathbf C)$ of $\mathfrak m(\mathbf C)$, it is enough to check
that the projection of $v_n$ into each simple root space of $(\mathfrak b
\cap \mathfrak m)(\mathbf C)$ is non-zero, and this follows from the
statement (proved above) that the image of $v$ under
$\mathfrak b(\mathbf C)
\hookrightarrow
\mathfrak p_\alpha(\mathbf C)
\twoheadrightarrow \mathfrak m_\alpha(\mathbf C)$ is regular in $\mathfrak
m_\alpha(\mathbf C)$ for every simple root $\alpha$ of $A$ in $M$. 

\subsection{End of the proof of the main theorem in case $A=T$}  We continue
with $u \in \mathfrak a(\mathcal O)$ and $x \in X^u$ as above. Let $M$ be
any Levi subgroup containing $A$. It remains to prove that if $x_P$ is
regular in $X_M^u$ for all $P \in \mathcal P(M)$ and if 
\begin{equation}\label{star} n(x,P,P')=n(u,P,P') 
\end{equation} for every adjacent pair $P,P' \in \mathcal P(M)$, then $x$ is
regular in $X^u$. We have already proved this in case $M=A$, and now we want
to reduce the general case to this special case. 

The equality \eqref{star} is equivalent to the equality
\begin{equation}\label{**}
\nu_M(x_P)-\nu_M(x_{P'})=n(u,P,P')\cdot \beta_{P,P'}. 
\end{equation} Fix $P \in \mathcal P(M)$ and sum \eqref{**} over the set of
neighboring pairs in a minimal gallery joining $P$ to its opposite $\bar P
\in \mathcal P(M)$. Doing this yields the equality 
\begin{equation}\label{***}
\nu_M(x_P)-\nu_M(x_{\bar P})=\sum_{\alpha \in R_N} \mathrm{val}_F(\alpha(u))
\cdot
\pi_M(\alpha^\vee), 
\end{equation}  where $\pi_M:X_*(A)\to\Lambda_M$ is the canonical surjection
and $R_N$ is the set of roots of $A$ in $\mathfrak n$. 

Fix a Borel subgroup $B_M$ in $M$ containing $A$ and let $B$ (respectively,
$B_1$) be the Borel subgroups in $\mathcal P(A)$ obtained as the inverse
image of $B_M$ under $P \to M$ (respectively,
$\bar P \to M$). Then \eqref{***} implies (see \ref{retractions}) that 
\begin{equation*}
\nu_A(x_B)-\nu_A(x_{B_1})\equiv \sum_{\alpha \in R_N}
\mathrm{val}_F(\alpha(u)) \cdot
\alpha^\vee  
\end{equation*} modulo the coroot lattice for $M$. Since $R_N$ is also the
set of roots that are positive on $B$ and negative on $B_1$, it follows that
\begin{equation*}
\nu_A(x_B)-\nu_A(x_{B_1})= \sum_{\alpha \in R_N} j_\alpha \cdot
\alpha^\vee  
\end{equation*} for some integers $j_\alpha$ such that $ 0 \le j_\alpha \le
\mathrm{val}_F(\alpha(u))$. (To prove this pick a minimal gallery joining
$B$ to $B_1$ and use the inequality stated in the main theorem for each
neighboring pair in the gallery.) Comparing this equality with the
congruence, we see that the linear combination 
\begin{equation}\label{4*}
 \sum_{\alpha \in R_N} (\mathrm{val}_F(\alpha(u))-j_\alpha) \cdot
\alpha^\vee  
\end{equation} maps to $0$ in $\Lambda_M$.

We get a basis for $\Lambda_M \otimes \mathbf R$ by taking the elements
$\beta_{P,P'}$ as $P'$ varies through the set of parabolic subgroups in
$\mathcal P(M)$ adjacent to $P$. Moreover for any
$\alpha \in R_N$ the image $\pi_M(\alpha^\vee)$ of $\alpha ^\vee$ in
$\Lambda_M$ is a non-negative linear combination of basis elements
$\beta_{P,P'}$ (with at least one non-zero coefficient). Therefore the fact
that \eqref{4*} maps to $0$ in $\Lambda_M$ means that
\begin{equation}\label{5*}
\nu_A(x_B)-\nu_A(x_{B_1})= \sum_{\alpha \in R_N} \mathrm{val}_F(\alpha(u))
\cdot
\alpha^\vee  
\end{equation}

By hypothesis $x_{\bar P}$ is regular. Therefore (transitivity of
retractions plus the part of our theorem we have already proved) for all
adjacent Borel subgroups $B_1,B_2 \in \mathcal P(A)$ such that $B_1,B_2
\subset \bar P$ we have
\begin{equation*}
\nu_A(x_{B_1})-\nu_A(x_{B_2})=  \mathrm{val}_F(\alpha_{B_1,B_2}(u)) \cdot
\alpha_{B_1,B_2}^\vee,  
\end{equation*} where $\alpha_{B_1,B_2}$ denotes the unique root that is
positive on $B_1$ and negative on $B_2$. Summing these equalities over
neighboring pairs in a minimal gallery joining $B_1$ to $\bar B$, we find
that
\begin{equation*}
\nu_A(x_{B_1})-\nu_A(x_{\bar B})= \sum_{\alpha \in R_M^+}
\mathrm{val}_F(\alpha(u)) \cdot
\alpha^\vee,  
\end{equation*} where $R_M^+$ denotes the set of roots of $A$ in $B_M$.
Adding this last equality to \eqref{5*}, we see that
\begin{equation}\label{6*}
\nu_A(x_{B})-\nu_A(x_{\bar B})= \sum_{\alpha \in R^+}
\mathrm{val}_F(\alpha(u)) \cdot
\alpha^\vee.  
\end{equation}

Now consider any minimal gallery $B=B_0,B_1,\dots,B_l=\bar B$ joining $B$ to
$\bar B$. Then
\begin{equation}\label{7*}
\nu_A(x_{B})-\nu_A(x_{\bar B})= \sum_{i=1}^l n(x,B_{i-1},B_i) \cdot
\alpha_i^\vee,  
\end{equation} where $\alpha_i$ is the unique root that is positive for
$B_{i-1}$ and negative for $B_i$. We know that $n(x,B_{i-1},B_i) \le
\mathrm{val}_F(\alpha_i(u))$ for all $i$. Subtracting \eqref{7*} from
\eqref{6*}, we find that $0$ is a non-negative linear combination of
positive roots; therefore each coefficient in this linear combination is
$0$, which means that
\begin{equation*}
 n(x,B_{i-1},B_i)=\mathrm{val}_F(\alpha_i(u)) 
  \end{equation*} for $i=1,\dots,l$.

Now consider any pair $B',B''$ of adjacent Borel subgroups in $\mathcal
P(A)$. After reversing the order of $B',B''$ if necessary we can find a
minimal gallery as above and an index $i$ such that 
$(B_{i-1},B_i)=(B',B'')$. Therefore
\begin{equation}\label{8*}
 n(x,B',B'')=\mathrm{val}_F(\alpha(u)), 
  \end{equation} where $\alpha$ is the unique root that is positive on $B'$
and negative on $B''$. Since both sides of \eqref{8*} remain unchanged when
$B',B''$ are switched, we see that \eqref{8*} holds for any adjacent pair
$B',B''$. By what we have already done, it follows that $x$ is regular in
$X^u$. 
\subsection{Proof of the main theorem in general} Now let  $M$ be any Levi
subgroup of $G$ containing $A$, and let $u$ be an integral regular
semisimple element of $\mathfrak g(F)$ that happens to lie in $\mathfrak
m(F)$. Let
$T=\mathrm{Cent}_{G_F}(u)$, a maximal torus in $M_F$. We choose a finite
extension $F'/F$ that splits $T$. 

We normalize the valuation $\mathrm{val}_{F'}$ on $F'$ so that uniformizing
elements in $F'$ have valuation $1$. Thus
$\mathrm{val}_{F'}(\epsilon)=[F':F]$. We write $X'$ for the set
$G(F')/G(\mathcal O_{F'})$. The inclusion $G(F) \hookrightarrow G(F')$
induces a canonical injection $X \hookrightarrow X'$.

For any $P \in \mathcal P(M)$ the diagram
$$\begin{CD} X @>r_P>> X_M \\ @VVV @VVV \\ X' @>r'_P>>X'_M 
\end{CD}$$ commutes, where the horizontal maps are retractions and the
vertical maps are the canonical injections. Moreover the diagram
\begin{equation*}
\begin{CD} X @>\nu_G>> \Lambda_G \\ @VVV @VVV \\ X' @>\nu'_G>>\Lambda_G 
\end{CD}
\end{equation*} commutes, where the left vertical map is the canonical
injection and the right vertical map is multiplication by $e:=[F':F]$.

For any $x \in X^u$ the image of $x$ in $X'$ lies in $(X')^u$, and $x$ is
regular in $X^u$ if and only if $x$ is regular in $(X')^u$. Indeed the
conjugacy class $\bar u_G(x)$ attached to $u$ and
$x$ is the same for $X$ and $X'$.

The torus $T$ is conjugate under $M(F')$ to $A$, so our theorem is true for
$T$ over $F'$. Therefore for $x \in X^u$ and adjacent $P=MN,P'=MN' \in
\mathcal P(M)$
\begin{equation}\label{2.6.1} e \cdot  n(x,P,P') \le \sum_{\alpha \in R_N
\cap R_{\bar N'}} \mathrm{val}_{F'}(\alpha(u)) \cdot m_\alpha, 
  \end{equation} and $x$ is regular in $X^u$ if and only if all of these
inequalities are equalities. (As before
$R_N$ denotes the set of roots of $A$ in $\mathfrak n$; the positive
integers $m_\alpha$ were defined in \ref{set-up}.) Dividing  by $e$, and
noting that the term indexed by
$\alpha$ depends only on the
$\Gamma$-orbit of $\alpha$, we find that \eqref{2.6.1} is equivalent to the
inequality
\begin{equation*}
  n(x,P,P') \le n(u,P,P').
\end{equation*} This completes the proof of  the theorem.

\providecommand{\bysame}{\leavevmode\hbox to3em{\hrulefill}\thinspace}


\begin{thebibliography}{Lan80}

\bibitem[Art76]{arthur76}
J.~Arthur, \emph{The characters of discrete series as orbital integrals},
  Invent. Math. \textbf{32} (1976), 205--261.

\bibitem[Bez96]{bezrukavnikov96}
R.~Bezrukavnikov, \emph{The dimension of the fixed point set on affine flag
  manifolds}, Math. Res. Lett. \textbf{3} (1996), 185--189.

\bibitem[KL88]{kazhdan-lusztig88}
D.~Kazhdan and G.~Lusztig, \emph{Fixed point varieties on affine flag
  manifolds}, Israel J. Math. \textbf{62} (1988), 129--168.

\bibitem[Lan80]{langlands80}
R.~P. Langlands, \emph{Base change for {$GL(2)$}}, Ann. of Math. Studies 96,
  Princeton University Press, 1980.

\end{thebibliography}
\end{document}